\documentclass[12pt]{article}

\usepackage{amssymb}

\begin{document}

\author{Antal Bege\\
"Babe\c s-Bolyai" University\\
Faculty of Mathematics\\
Str. Kog\^alniceanu, Nr. 1\\
3400 Cluj-Napoca, Romania}
\title{\bf A generalization of Apostol's M\"obius functions of order $k$}
\date{September 12, 2000}
\maketitle

\begin{abstract}
Apostol's M\"obius functions $\mu_k(n)$ of order $k$ are generalized to depend 
on a~second integer parameter $m\geq k$. Asymptotic formulas are obtained for 
the partial sums of these generalized functions.
\end{abstract}

\section{Introduction}
M\"obius functions of order $k$, introduced by T. M. Apostol \cite{apostol1}, 
are defined by the formulas
$$
\mu_k(n)=
\left\{
\begin{array}{cl}
1& \mbox{ if }n=1,\\
0& \mbox{ if }p^{k+1}\mid n\mbox{ for some prime } p,\\
(-1)^r& \mbox{ if } n=p_1^k\cdots p_r^k\prod\limits_{i>r}p_i^{\alpha_i},
\quad \mbox{ with } 0\leq \alpha_i<k,\\
1& \mbox{ otherwise. }  
\end{array}
\right.
$$
In \cite{apostol1} Apostol obtained the asymptotic formula
\begin{equation}
\label{11}
\sum_{n\leq x}\mu_k(n)=A_kx+O(x^{\frac{1}{k}}\log x),
\end{equation}
where
$$
A_k=\prod_p\left(1-\frac{2}{p^k}+\frac{1}{p^{k+1}}\right).
$$
Later, Suryanarayana \cite{suryanarayana1} showed that, on the assumption of 
the Riemann hypothesis, the error term in (\ref{11}) can be improved to
\begin{equation}
\label{12}
O\left(x^{\frac{4k}{4k^2+1}}\omega (x)\right),
\end{equation}
Where
$$
\omega (x)=\mbox{exp}\{A\log x (\log \log x)^{-1}\}
$$
for some positive constant $k$.
\\
\\
This paper generalizes M\"obius functions of order $k$ and establishes 
asymptotic formulas for their partial sums.

\section{Preliminary lemmas}

The generalization in question is denoted by $\mu_{k,m}(n)$, where $1<k\leq m$.
\\
If $m=k$, $\mu_{k,k}(n)$ is defined to be $\mu_k(n)$, and if $m>k$ the 
function is defined as follows:
\begin{equation}
\label{21}
\mu_{k,m}(n)=
\left\{
\begin{array}{cl}
1& \mbox{ if } n=1,\\
1& \mbox{ if } p^{k}\nmid n\mbox{ for each prime } p,\\
(-1)^r& \mbox{ if } n=p_1^m\cdots p_r^m\prod\limits_{i>r}p_i^{\alpha_i},
\quad \mbox{ with } 0\leq \alpha_i<k,\\
0& \mbox{ otherwise. }  
\end{array}
\right.
\end{equation}
This generalization, like Apostol's $\mu_k(n)$, is a~multiplicative function 
of $n$, so it is determined by its values at the prime powers. We have
$$
\mu_k(p^{\alpha })=\left\{
\begin{array}{rcl}
1& \mbox{ if }& 0\leq \alpha <k,\\
-1& \mbox{ if }& \alpha =k,\\
0& \mbox{ if }& \alpha >k,
\end{array}
\right.
$$
whereas
\begin{equation}
\label{22}
\mu_{k,m}(p^{\alpha })=\left\{
\begin{array}{rcl}
1& \mbox{ if }& 0\leq \alpha <k,\\
0& \mbox{ if }& k\leq \alpha <m,\\
-1& \mbox{ if }& \alpha =m,\\
0& \mbox{ if }& \alpha >m,
\end{array}
\right.
\end{equation}

\noindent
{\bf Lemma 2.1}
{\it
For $k\leq m$ we have
\begin{equation}
\label{23}
\mu_{k,m}(n)=\sum_{
\begin{array}{l}
\delta d^m=n\\
(d,\delta )=1
\end{array}}
\mu (d)\; q_k(\delta )
\end{equation}
where $\mu (n)$ is the M\"obius function and $q_k(n)$ is the caracteristic 
function of the k-free integers:}
$$
q_k(n)=\left\{
\begin{array}{lcl}
0& \mbox{ if }& p^k|n \mbox{ for some prime }p\\
1& \mbox{ if }& p^{\alpha }|n \mbox{ implies }\alpha <k
\end{array}
\right. .
$$

\noindent
{\bf Proof}
\\
Because $\mu (n)$ and $q_k(n)$ are multiplicative functions of $n$, the sum in 
the lemma is also multiplicative, so to complete the proof we simply note that 
when $n=p^{\alpha }$ the sum has the values indicated in (\ref{22}).
\\

\noindent
The next two lemmas, proved in \cite{suryanarayana2}, involve the following 
functions:
$$
\theta(n)=\mbox{ the number of square-free divisors of }n,
$$
$$
\psi_k(n)=n\prod_{p|n}\left(1+\frac{1}{p}+\cdots +\frac{1}{p^{k-1}}\right),
$$
where $k$ is an integer $\geq 2$,
$$
\delta _k(x)=\mbox{ exp }\{-A\; k^{-\frac{8}{5}}\log^{\frac{3}{5}}x\; 
(\log \log x)^{-\frac{1}{5}}\},
$$
where $A>0$ is an absolute constant,
$$
\omega_k(x)=\mbox{ exp }\{B_k\; \log x\; (\log \log x)^{-1}\},
$$
where $B_k$ is a positive constant.
\\

\noindent
{\bf Lemma 2.2}
{\it
For $x\geq 3$ we have
\begin{equation}
\label{24}
Q_k(x,n)=\sum_{
\begin{array}{c}
r\leq x\\
(r,n)=1
\end{array}}
q_k(r)=
\frac{xn}{\zeta (k) \psi_k(n)}+
0\left(\theta (n)x^{\frac{1}{k}}\delta_k(x)\right)
\end{equation}
uniformly in $x$, $n$ and $k$.}
\\

\noindent
{\bf Lemma 2.3}
{\it
If the Riemann hypothesis is true, then for $x\geq 3$ we have
\begin{equation}
\label{25}
Q_k(x,n)=\sum_{
\begin{array}{c}
r\leq x\\
(r,n)=1
\end{array}}
q_k(r)=
\frac{xn}{\zeta (k) \psi_k(n)}+
0\left(\theta (n)x^{\frac{2}{2k+1}}\omega_k(x)\right)
\end{equation}
uniformly in $x$, $n$ and $k$.}
\\

\noindent
Our derivation of an asymptotic formula for the summatory function of 
$\mu_{k,m}(n)$ will also make use of the following lemma.
\\

\noindent
{\bf Lemma 2.4}
{\it
For $k\geq 2$ we have
\begin{equation}
\label{26}
\sum_{d|n}\frac{\mu(d)\psi_{k-1}(d)}{d\psi_k(d)}=
\frac{n}{\psi_k(n)}.
\end{equation}}

\noindent
{\bf Proof}
\\
Both sides of (\ref{26}) are multiplicative functions of $n$ so it suffices 
to verify the equation when $n$ is a~prime power. If $n=p^{\alpha }$ we have
$$
\sum_{d|p^{\alpha }}\frac{\mu(d)\psi_{k-1}(d)}{d\psi_k(d)}=
1-\frac{1+\frac{1}{p}+\cdots +\frac{1}{p^{k-2}}}
{p\left(1+\frac{1}{p}+\cdots +\frac{1}{p^{k-1}}\right)}=
$$
$$
=\frac{p^{\alpha }}
{p^{\alpha }\left(1+\frac{1}{p}+\cdots +\frac{1}{p^{k-1}}\right)}=
\frac{n}{\psi_k(n)}.
$$

\noindent
The next lemma is proved in \cite{suryanarayana3}.
\\
\\
{\bf Lemma 2.5}
{\it 
For $x\geq 3$, $n\geq 1$, and every $\epsilon >0$ we have}
\begin{equation}
\label{27}
L_n(x)\equiv \sum_{
\begin{array}{c}
r\leq x\\
(r,n)=1
\end{array}}
\frac{\mu(r)}{r}=
O\left(\sigma ^*_{-1+\epsilon}(n)\delta (x)\right)
\end{equation}
{\it uniformly in $x$ and $n$, where $\sigma ^*_{\alpha }(n)$ is the sum of the 
$\alpha $th powers of the square-free divisors of $n$, and
$$
\delta (x)=\mbox{ exp }\{-A\; \log^{\frac{3}{5}}x\; 
(\log \log x)^{-\frac{1}{5}}\},
$$
for some absolute constant $A>0$.}
\\
\\
We note that if $\alpha <0$ we have $\sigma^*_{\alpha }(n)\leq \sigma^*_0(n)=
\theta (n)$. Also, $x^{\epsilon }\delta (x)$ is an increasing function of $x$ 
for every $\epsilon >0$ and $x>x_0(\epsilon )$. Using the method described in 
\cite{suryanarayana3}, 
it can be shown that if the Riemann hypothesis is true the factor 
$\delta (x)$ in the error term in (\ref{27}) can be replaced by 
$\omega (x)x^{-\frac{1}{2}}$, where
$$
\omega (x)=\mbox{ exp }\{A\; \log x\; (\log \log x)^{-1}\},
$$
for some absolute constant $A>0$.
\\
\\
{\bf Lemma 2.6}
{\it
For $x\geq 3$ and every $\epsilon >0$ we have
\begin{equation}
\label{28}
\sum_{
\begin{array}{c}
r\leq x\\
(r,n)=1
\end{array}}
\frac{\mu(r)}{\psi_k(r)}=
0\left(\sigma^*_{-1+\epsilon}(n)\delta(x)\right).
\end{equation}
uniformly in $x$, $n$ and $k$}.
\\
\\
{\bf Proof}
We write
$$
\sum_{
\begin{array}{c}
r\leq x\\
(r,n)=1
\end{array}}
\frac{\mu(r)}{\psi_k(r)}=
\sum_{
\begin{array}{c}
r\leq x\\
(r,n)=1
\end{array}}
\frac{\mu(r)}{r}\frac{r}{\psi_k(r)},
$$
then use (\ref{26}) to obtain
$$
\sum_{
\begin{array}{c}
r\leq x\\
(r,n)=1
\end{array}}
\frac{\mu(r)}{\psi_k(r)}=
\sum_{
\begin{array}{c}
d\delta\leq x\\
(d,\delta )=1\\
(d\delta ,n)=1
\end{array}}
\frac{\mu^2(d)\mu (\delta )\psi_{k-1}(d)}{d^2\delta \psi_k(d)}=
$$
$$
=\sum_{
\begin{array}{c}
d\leq x\\
(d,n)=1
\end{array}}
\frac{\mu^2(d)\psi_{k-1}(d)}{d^2\psi_k(d)}\; 
\sum_{
\begin{array}{c}
\delta\leq \frac{x}{d}\\
(\delta ,nd)=1
\end{array}}
\frac{\mu(\delta )}{\delta}=
$$
$$
=\sum_{
\begin{array}{c}
d\leq x\\
(d,n)=1
\end{array}}
\frac{\mu^2(d)\psi_{k-1}(d)}{d^2\psi_k(d)}\; 
L_{dn}\left(\frac{x}{d}\right).
$$
Using (\ref{27}) and the inequality $\psi_{k-1}(d)\leq \psi_k(d)$ we find that 
the last sum is
$$
O(\sigma^*_{-1+\epsilon}(n)\; 
\sum_{
\begin{array}{c}
d\leq x\\
(d,n)=1
\end{array}}
\frac{\mu^2(d)\delta\left(\frac{x}{d}\right)
\sigma^*_{-1+\epsilon }(d)}{d^2}).
$$ 
Because $x^{\epsilon '}\delta (x)$ increases for every $\epsilon '>0$, 
we have $\left(\frac{x}{d}\right)^{\epsilon '}\delta\left(\frac{x}{d}\right) 
\leq 
x^{\epsilon '}\delta (x)$ so $\delta \left(\frac{x}{d}\right)\leq d^{\epsilon '}\delta (x)$
and the foregoing $0$-term is
$$
O (\sigma^*_{-1+\epsilon}(n)\delta(x)\; 
\sum_{
\begin{array}{c}
d\leq x\\
(d,n)=1
\end{array}}
\frac{\mu^2(d)\sigma^*_{-1+\epsilon }(d)}{d^{2-\epsilon'}}).
$$ 
But $\sigma^*_{-1+\epsilon}(d)\leq \tau(d)=O(d^{\epsilon '})$ for every
$\epsilon '>0$. If we choose $\epsilon '<\frac{1}{2}$ the last sum is
$O(1)$ and we obtein (\ref{28}).
\\
\\
Applying (\ref{28}) together with the remark following Lemma 2.5 we obtain
\\
\\
{\bf Lemma 2.7}
{\it 
If the Riemann hypothesis is true, then for $x\geq 3$, $n\geq 1$ and every 
$\epsilon >0$ we have
\begin{equation}
\label{29}
\sum_{
\begin{array}{c}
r\leq x\\
(r,n)=1
\end{array}}
\frac{\mu(r)}{\psi_k(r)}=
O\left(\sigma^*_{-1+\epsilon}(n)\omega(x)x^{-\frac{1}{2}}\right).
\end{equation}
uniformly in $x$, $n$ and $k$}.
\\
\\
Applying partial summation in (\ref{28}) we obtain
\\
\\
{\bf Lemma 2.8}
{\it
For $x\geq 3$, $k\geq 2$, and every $\epsilon >0$ we have
\begin{equation}
\label{210}
\sum_{
\begin{array}{c}
r\leq x\\
(r,n)=1
\end{array}}
\frac{\mu(r)}{\psi_k(r)r^{k-1}}=
O\left(\sigma^*_{-1+\epsilon}(n)\delta(x)x^{1-k}\right).
\end{equation}
uniformly in $x$, $n$ and $k$}.
\\
\\
{\bf Note}: If the Riemann hypothesis is true, the error term in Lemma 2.8 
holds with $\delta (x)x^{1-k}$ replaced by $\omega (x) x^{\frac{1}{2}-k}$.

\section{Main results}

{\bf Theorem 3.1}
{\it
For $x\geq 3$ and $m>k\geq 2$ we have
\begin{equation}
\label{31}
\sum_{
\begin{array}{c}
r\leq x\\
(r,n)=1
\end{array}}
\mu_{k,m}(r)=
\frac{xn^2\; \alpha_{k,m}}{\zeta(k)\psi_k(n)\alpha_{k,m}(n)}+
0\left(\theta(n)x^\frac{1}{k}\delta(x)\right).
\end{equation}
uniformly in $x$, $n$ and $k$, where
$$
\alpha_{k,m}=\prod_p\left(1-\frac{1}{p^{m-k+1}+p^{m-k+2}+\cdots +p^m}\right)
$$
and
$$
\alpha_{k,m}(n)=n\prod_{p|n}\left(1-
\frac{1}{p^{m-k+1}+p^{m-k+2}+\cdots +p^m}\right).
$$}
\\
\\
{\bf Proof}
By (\ref{23}) and (\ref{24}) we have
$$
\sum_{
\begin{array}{c}
r\leq x\\
(r,n)=1
\end{array}}
\mu_{k,m}(n)=
\sum_{
\begin{array}{l}
\delta d^m\leq x\\
(d,\delta )=1\\
(d\delta ,n)=1
\end{array}}
\mu (d)\; q_k(\delta )=
$$
$$
=
\sum_{
\begin{array}{l}
d\leq x^{\frac{1}{m}}\\
(d,n)=1
\end{array}}
\mu (d)
\sum_{
\begin{array}{l}
\delta \leq \frac{x}{d^m}\\
(\delta ,dn)=1
\end{array}}
q_k(\delta )=
$$
$$
=
\sum_{
\begin{array}{l}
d\leq x^{\frac{1}{m}}\\
(d,n)=1
\end{array}}
\mu (d)
Q_k\left( \frac{x}{d^m}, dn\right)=
$$
$$
=\sum_{
\begin{array}{l}
d\leq x^{\frac{1}{m}}\\
(d,n)=1
\end{array}}
\mu (d)\left\{
\frac{\left(\frac{x}{d^m}\right)dn}{\zeta (k) \psi_k(dn)}+
0\left(\theta (dn)\frac{x^{\frac{1}{k}}}{d^{\frac{m}{k}}}
\delta \left(\frac{x}{d^m}\right)
\right)
\right\}=
$$
$$
=
\frac{xn}{\zeta(k)\psi_k(n)}
\sum_{
\begin{array}{c}
d=1\\
(d,n)=1
\end{array}}
^{\infty }
\frac{\mu(d)}{d^{m-1}\psi_k(d)}-
\frac{xn}{\zeta (k)}{\psi_k(n)}
\sum_{
\begin{array}{c}
d>x^{\frac{1}{m}}\\
(d,n)=1
\end{array}}
\frac{\mu (d)}{d^{m-1}\psi_k(d)}+
$$
$$
+
O(
\theta (n)x^{\frac{1}{k}-\epsilon}
\sum_{
\begin{array}{c}
d\leq x^{\frac{1}{m}}\\
(d,n)=1
\end{array}}
\frac{\delta\left(\frac{x}{d^k}\right)\left(\frac{x}{d^k}\right)\mu^2(d)
\theta(d)}
{d^{\frac{m}{k}-\epsilon 'k}}
).
$$
Using the Euler product representation for absolutely convergent series of 
multiplicative terms \cite{apostol2} we have
$$
\sum_{
\begin{array}{c}
d=1\\
(d,n)=1
\end{array}}
^{\infty }
\frac{\mu(d)}{d^{m-1}\psi_k(d)}=
\prod_{p\mid n}\left(1-
\frac{1}{p^{m-k+1}+\cdots +p^m}\right)=
\frac{\alpha_{k,m}}{\alpha_{k,m}(n)}.
$$
Now use (\ref{210}) and the fact that $\delta(x)x^{\epsilon '}$ is increasing
for all $\epsilon '>0$, then choose $\epsilon >0$ so that 
$\frac{m}{k}-\epsilon 'k>1+\epsilon $ and we obtain (\ref{31}).
\\
When $n=1$, Theorem 3.1 gives the following corollary for $x\geq 3$ and 
$m>k\geq 2$:
\begin{equation}
\label{32}
\sum_{r\leq x}
\mu_{k,m}(r)=
\frac{x}{\zeta(k)}\alpha_{k,m}+
0\left(x^{\frac{1}{k}}\delta(x)\right)
\end{equation}
uniformly in $x$ and $k$.
\\
Applying the method used to prove Theorem 1, and making use of (\ref{25}) and 
Lemma 2.9 we get
\\
\\
{\bf Theorem 3.2}
\\
{\it If the Riemann hypothesis is true, then for $x\geq 3$ and $m>k\geq 2$ 
we have
\begin{equation}
\label{33}
\sum_{
\begin{array}{c}
r\leq x\\
(r,n)=1
\end{array}}
\mu_{k,m}(r)=
\frac{xn^2\; \alpha_{k,m}}{\zeta(k)\psi_k(n)\alpha_{k,m}(n)}+
0\left(\theta(n)x^{\frac{2}{2k+1}}\omega(x)\right).
\end{equation}
uniformly in $x$, $n$ and $k$.}
\\
\\
In particular, if n=1 we have
\begin{equation}
\label{34}
\sum_{r\leq x}
\mu_{k,m}(r)=
\frac{x}{\zeta(k)}\alpha_{k,m}+
0\left(x^{\frac{2}{2k+1}}\omega(x)\right)
\end{equation}
uniformly in $x$ and $k$.

\section{Conjectures}

Suryanarayana raised the question of improving the error term in Apostol's 
asymptotic formula (\ref{11}), and notes that no improvement seems possible 
by this method. Our method gives no improvement in the error term but it does suggest the 
following conjectures:
\\
For $x\geq 3$, $n\geq 1$ and $k\geq 2$ 
we have
\begin{equation}
\label{41}
\sum_{
\begin{array}{c}
r\leq x\\
(r,n)=1
\end{array}}
\mu_{k}(r)=
\frac{xn^2\; \alpha_{k,k}}{\zeta(k)\psi_k(n)\alpha_{k,k}(n)}+
0\left(\theta(n)x^{\frac{1}{k}}\delta(x)\right).
\end{equation}
uniformly in $x$, $n$ and $k$.
\\
In particular, when $n=1$ the conjecture is
\begin{equation}
\label{42}
\sum_{r\leq x}
\mu_{k}(r)=
\frac{x}{\zeta(k)}\alpha_{k,k}+
0\left(x^{\frac{1}{k}}\delta(x)\right)
\end{equation}
uniformly in $x$ and $k$.
\\
If the Riemann hypothesis is true, the conjectured formulas are
\begin{equation}
\label{43}
\sum_{
\begin{array}{c}
r\leq x\\
(r,n)=1
\end{array}}
\mu_{k}(r)=
\frac{xn^2\; \alpha_{k,k}}{\zeta(k)\psi_k(n)\alpha_{k,k}(n)}+
0\left(\theta(n)x^{\frac{2}{2k+1}}\omega(x)\right).
\end{equation}
uniformly in $x$, $n$ and $k$, for $x\geq 3$, $n\geq 1$ and $k\geq 2$.
\\
In particular, when $n=1$ the conjecture is
\begin{equation}
\label{44}
\sum_{r\leq x}
\mu_{k}(r)=
\frac{x}{\zeta(k)}\alpha_{k,k}+
0\left(x^{\frac{2}{2k+1}}\omega(x)\right)
\end{equation}
uniformly in $x$ and $k$.
\\
\\
It should be noted that $\alpha_{k,k}=\zeta (k)A_k$, where $A_k$ is Apostol's 
constant in (\ref{11}), so the leading term in (\ref{42}) and (\ref{44}) is 
the same as that in (\ref{11}).

\begin{center}
ACKNOWLEDGEMENT
\end{center}

The author wishes to thank Professor Tom M. Apostol for his help in preparing 
this paper for publication.


\begin{thebibliography}{99}
\bibitem{apostol1}
T. M. Apostol, 
M\"obius functions of order $k$, 
{\it Pacific Journal of Math.}, 
{\bf 32} (1970), 21-17.
\bibitem{apostol2}
T. M. Apostol, 
{\it Introduction to Analytic Number Theory}, 
Undergraduete texts in Mathematics, Springer Verlag, New-York, 1976.
\bibitem{suryanarayana1}
D. Suryanarayana, 
On a theorem of Apostol concerning Mobius functions of order $k$, 
{\it Pacific Journal of Math.}, 
{\bf 68} (1977), 277-281.
\bibitem{suryanarayana2}
D. Suryanarayana, 
Some more remarks on uniform O-estimates for k-free integers, 
{\it Indian J. Pure Appl. Math.}, 
{\bf 12} (11) (1981), 1420-1424.
\bibitem{suryanarayana3}
D. Suryanarayana and P. Subrahmanyam, 
The maximal k-free divisor of m which is prime to n, 
{\it Acta Math. Acad. Sci. Hung.}, 
{\bf 33} (1979), 239-260.
\end{thebibliography}
\end{document}